\documentclass{article}
\usepackage[cp1251]{inputenc}    % Перешли на кодировку Windows!!!
\usepackage[english, russian]{babel}    % Переносы через Babel (обязательно!)

\usepackage{amsmath}            % Если используются возможности пакета
\usepackage{amsfonts,amssymb}   % Если есть необходимость выбора
\usepackage{cite}
\newcommand{\doc}{{\it{Доказательство.}}\ }

\newtheorem{sled}{Следствие}
\newtheorem{prop}{Предложение}
\newtheorem{zam}{\mbox{З а м е ч а н и е}\ \ \\ }

\newcommand{\rav}{\stackrel{\triangle}{=}}
\newcommand{\ph}{\varphi}

\newcommand{\epsi}{\varepsilon}

\newcommand{\rref}[1]{$(\ref{#1})$}
\newcommand{\mm}[1]{{\bf{#1}}}
\newcommand{\ct}[1]{{\mathcal{#1}}}
\newcommand{\fr}[1]{{\mathfrak{#1}}}
\newcommand{\td}[1]{\tilde{#1}}
\newcommand{\bo}{{\hfill {$\Box$}}}
\usepackage[active]{srcltx}
\begin{document}
 \centerline{\Large Задачи управления на бесконечном промежутке и}
 \centerline{ \Large  устойчивость сопряженной переменной}
 \vskip+0,5cm

%\title{Принцип максимума для задач управления на бесконечном промежутке}{Хлопин~Д.В.%
%\protect\footnote{Работа поддержана грантом РФФИ №~09-01-00436 и
% программой президиума РАН <<Математическая теория управления>>.}}%
%[Хлопин~Д.В. Задачи  на бесконечном промежутке]%of the adjoint

%An in?nite-horizon optimal control problem and the stability of the adjoint variable (in Russian).

%In an in?nite-horizon optimal control problem the stability of the adjoint variable  
% implying the vanishing of the adjoint variable at in?nity along optimal solution.

%% Необязательный аргумент - для включения в содержание сборника.
%% Применяется в случае наличия в заголовке принудительных переносов, сносок и т.п.
%\vspace{-1.75\baselineskip}

\centerline{ Хлопин Д.В. (ИММ УрО РАН, Екатеринбург)}
\vspace{\baselineskip}

\centerline{\it e-mail:\ khlopin@imm.uran.ru } \vspace{\baselineskip}

  {\bf Определения и обозначения.}
  Пусть задано
 метрическое пространство $\mm{X}\rav\mm{R}^m$.
 Определим $\mm{T}\rav\{t\in\mm{R}\,|\,t\geq 0\}.$
 Всюду на пространствах функций на $\mm{T}$
  рассматривается компактно-открытая топология.
 В частности, введем $\fr{X}\rav C(\mm{T},\mm{X}).$
 Обозначим также через $\Omega$ семейство тех функций $\omega\in C(\mm{T},\mm{R}),$
 для которых $\lim_{t\to\infty}\omega(t)=0.$
% На $C(\mm{T},\mm{X})$ введем компактно-открытую топологию.

% Пусть также задано
% некоторое множество $\mm{K}\subset C(\mm{T},\mm{X}),$ причем для
% всех $t\in\mm{T}$ множества $\mm{K}_t\rav{x_{[0,t]}\,|\,x\in\mm{K}}$
% компактны.
 %Через $exp\,\mm{K}$ обозначим всевозможные непустые
 %множества подмножества $A$ множества $\mm{K}$ со свойством: для всех
 %$t\in\mm{T}$ $A_t\rav\{x_{}[0,t]\,|\,x\in A\}$ замкнуто. Пусть дано
 %некоторое множество $\mm{V}\subset exp\,\mm{K}.$
% Это множество может играть роль ограничений на траектории.

  {\bf % Исходная
  Управляемая система.}
   Пусть дана управляемая система
   \begin{equation}
   \label{sys}
    \dot{x}=f(t,x,u),\ x(0)=0,\ t\in\mm{T}\rav\mm{R}_{\geq 0},\ u\in P\in (comp)(\mm{R}^p).
   \end{equation}

 Условие $\bf{I}:$ Для всякого $t\in\mm{T}$ для всякой измеримой функции $u\in B([0,t],P)$ найдется
 единственное решение $\ph_t(u)\in C([0,t],\mm{X})$. Более того, эта
 зависимость $u\mapsto\ph_t(u)$ непрерывна для всех $n\in\mm{N}$.

 Заметим, что это условие фактически гарантирует как единственность
 траектории, так и ее неуход на бесконечность за конечное время.
 С другой стороны для обеспечения этого условия достаточно
 принять:

 Условие $\bf{I}a:$ Условия Каратеодори, то есть:\
 1) для произвольных $(x,u)\in \mm{X}\times P$
 функция $(f(t,x,u)|t\in \mm{T})$
 измерима;\
 2) для каждого $t\in \mm{T}$ функция
 $(f(t,x,u)|(x,u)\in \mm{X} \times P)$
 непрерывна;

 Условие $\bf{I}b:$
 Для всякого $n\in\mm{N}$ найдется
 такая суммируемая функция
 $L:[0,n]\mapsto \mm{T}$, что функция
 $(f(t,x,u)|x\in \mm{X})$
 липшицева по $x$ с константой $L(t)$ для почти всех $t\in [0,n]$ и всех $u\in
 P$;

 %Условие локальной липшицевости  по $\max t,x,u$.

 Условие $\bf{I}c:$ условие продолжимости всех
 решений на $\mm{T}$, например
 условие подлинейного (по $x$) роста.
%  для некоторой суммируемой на всяком компакте функции  $l\in B(\mm{T},\mm{T})$
%   и непрерывной, но несуммируемой на $\mm{T}$ функции
%    $R:\mm{T}\to\mm{R}_{> 0}$ выполнено
%  $$  ||f(t,x,u)||\leq l(t)/R(||x||)\ \ \ \forall (t,x,u)\in \mm{T}\times\mm{X}\times P;$$

% Функции $\bf{I}c$ может выполнять  условие подлинейного
% роста.

  Введем множество $\fr{U}\rav B(\mm{T},P).$
 В силу условия $\bf{I}$ отображение,   %непрерывную функцию,
 сопоставляющее
 всякому  $u\in\fr{U}$ его решение $\ph(u)\in\fr{X}$,
 непрерывно.

% Примем $\Phi_n\rav\ph_n^{B([0,n],P)}.$

 {\bf Задача управления.} Пусть поставлена задача максимизации  на траекториях \rref{sys}
 функционала
 \begin{equation}
   \label{opt}
   J(u)\rav\int_{\mm{T}} g(t,\ph(u)(t),u(t))dt\to\max,
   \end{equation}
 Для существования и непрерывности функционала $J$ на допустимых парах $(u,\ph(u))$
 предположим:

 Условие $\bf{II}:$ для $g$ выполнены $\bf{I}a)-c)$, а кроме того
 для некоторой функции $\omega\in\Omega$
 при любых $u\in \fr{U}$ выполнено
 $$\int_{T}^{\infty}|g(t,\ph(x)(t),u(t))|dt \leq \omega(T),\ \ \forall T\in\mm{T}.$$

 В частности теперь $J$ ограничено сверху.

 Определим функцию Гамильтона-Понтрягина
  $\ct{H}:\fr{X}\times\mm{T}\times\fr{U}\times\mm{R}\times\fr{X}\mapsto\mm{R}$
  правилом: $\ct{H}(x,t,u,\lambda,\psi)\rav\psi'f(t,x,u(t))+\lambda g(t,x,u(t)).$
 Введем соотношения
\begin{equation}
   \label{sys_x}
       \dot{x}(t)= f(t,x(t),u(t));
   \end{equation}
\begin{equation}
   \label{sys_psi}
       \dot{\psi}(t)\in-{\partial_x \ct{H}(x(t),t,u(t),\lambda,\psi(t))},
       \ \ ||\psi(0)||^2+\lambda^2=1;
   \end{equation}
\begin{equation}
   \label{maxH}
        \ct{H}(x(t),t,u(t),\lambda,\psi(t))=\sup_{p\in
        P}\ct{H}(x(t),t,p,\lambda,\psi(t)).
   \end{equation}

 Заметим, что в силу свойств {\bf{I,II}}
 включение \rref{sys_psi} полунепрерывно сверху, имеет измеримую
 мажоранту, а следовательно в силу  \cite[Теорема 4.1]{tovst}
 всякое решение этих соотношений продолжимо на все $\mm{T}.$

{\bf Существование.}
 Покажем в
 условиях $\bf{I}-\bf{II}$ существование оптимального решения,
 расширив множество управлений.

 Пусть {$\Pi([0,n], P)$}
  --- оснащенное топологией *-слабой
  сходимости множество всех слабо измеримых отображений из $[0,n]$
  в множество вероятностных мер Радона над $P.$
  Определим топологическое пространство $\Pi(\mm{T}, P)$ как обратный предел (
  \cite[III.1.30]{phillvv}, \cite[2.5.6]{en})%\cite{phillvv})
   $$\Pi(\mm{T}, P)=\lim_{\leftarrow}\mm{S}=\lim_{\leftarrow}\{
  \Pi([0,n], P),\pi^{[0,n]}_{[0,m]},\mm{N}\}.$$
   (здесь отрезки направлены отношением $\subset$ и
   $\pi^{K'}_{K''}(\eta)=\eta|_{K''}$ для всех $\eta\in \Pi(K', P)$).
  В частности элементы $\Pi(\mm{T}, P)$ ---
   такие функции  $\eta$,
  что  $\eta|_{[0,n]}\in \Pi([0,n], P)$ для всякого  $n\in\mm{N}$, а
  множество $\ct{A}\subset\Pi(\mm{T}, P)$ замкнуто, только если
  для всех $t\in\mm{T}$ образ $\pi_{[0.t]}(\ct{A})$ замкнут в
  $\Pi([0,t], P)$.

  Каждой  $\eta\in {\Pi(\mm{T}, P)}$ можно сопоставить
  $\td\ph(\eta)\in C(\mm{T},\mm{X})$ как
  решение уравнения
   $$\dot{x}=\int_{P} f(\tau,x(\tau),u)\,\eta(t)(du),\ x(0)=x_0$$
  и функционал $\td{J}(\eta)=\int_{\mm{T}}\int_{P} g(\tau,\td\ph(\eta)(\tau),u)\,\eta(t)(du)\,dt.$

  Как обратный предел компактов $\Pi([0,n], P)$ само $\Pi(\mm{T}, P)$ также компакт
  (Теорема Куроша \cite[III.1.13]{phillvv}). Также как $B([0,n],P)$
  всюду плотно вкладывалось в $\Pi([0,n], P)$, теперь и $\fr{U}$ всюду плотно
   вкладывается в $\Pi(\mm{T}, P)$ (\cite[III.1.27]{phillvv}), а  отображения
   $\td\ph$ и $\td{J}$ являются фактически продолжениями по непрерывности
  функционалов $\ph$ и $J$ с всюду плотного множества $\fr{U}$ на компакт
  $\Pi(\mm{T}, P)$. Отсюда обобщенная задача
 максимизировать
 функционал $\td{J}(\eta)$
  является релаксацией исходной, то есть
% непрерывно продолжается на
%    $\Pi(\mm{T}, P)$,
 % можно непрерывно продолжить, при этом
 $\sup_{u\in\fr{U}}J(u)=\max_{\eta\in\Pi(\mm{T}, P)}\td{J}(\eta)$,
 у
  обобщенного аналога \rref{opt} в
 условиях $\bf{I}-\bf{II}$ найдётся оптимальное решение $\mu_0\in\Pi(\mm{T}, P),$
 всякому такому управлению можно сопоставить последовательность
 управлений $(u_n)_{n\in\mm{N}}\in\fr{U}^\mm{N},$
 сходящуюся к $\mu_0$ в $\Pi(\mm{T}, P)$, а у всякой реализующей супремум $J$ последовательности
 $(u_n)_{n\in\mm{N}}\in\fr{U}^\mm{N}$ найдется предельная точка, на которой реализуется максимум
 $\td{J}$.

 Перейдя теперь к управлениям
  Гамкрелидзе (см. \cite{ga,kr_as}) можно обеспечить

 Условие $\bf{III}$ (\cite[(A2)]{kr_as}).  Для всякого $(t,x)\in\mm{T}\times \mm{X}$
выпукло  множество
  $\{(z',f(t,x,u))\in\mm{R}^{m+1}\,|\,z\in \langle -\infty,g(t,x,u)],u\in P\}.$
% поскольку в этом случае $\td{\ph}(\Pi(\mm{T}, P))=\ph(\Pi(\mm{T}, P))$
% и для всякой $\eta\in\Pi(\mm{T}, P)$ найдется $u\in B(\mm{T},P)$
% со свойством $\td\ph(u)=\ph(\eta)$, при этом $g(t,\td\ph(u)(t),u(t))\leq$

% Можно привести условие, при котором $\mu_0$ фактически элемент $\fr{U},$
  %это обобщенное оптимальное решение лежит в
%  $\fr{U},$
% более того,
  Тогда будут замкнутыми образы $J_\tau(\fr{U}))$ для всякого $\tau\in\mm{T}$
  (здесь $J_\tau(u)\rav\int_{[0,\tau]}g(t,\ph(u)(t),u)dt$),
  откуда следует компактность $J(\fr{U})),$ тем самым
   показано

\begin{sled}[\cite{bald}]
    В условиях $\bf{I}-\bf{III}$
    для задачи \rref{opt} всегда существует оптимальное
    управление
    $u\in\fr{U}\rav B(\mm{T},P).$
\end{sled}

   Аналогично можно ввести систему для обобщенных решений.
\begin{equation}
   \label{sys_x_}
       \dot{x}_*(t)= \int_{P} f(x_*(t),u) \eta(t)(du);
   \end{equation}
\begin{equation}
   \label{sys_psi_}
       \dot{\psi}(t)\in-\int_{P} {\partial_x \ct{H}(x_*(t),t,u_*(t),\lambda,\psi(t))}d\eta(t)(du),
       \ \ ||\psi(0)||^2+\lambda^2=1;
   \end{equation}
\begin{equation}
   \label{maxH_}
        \int_P\ct{H}(x_*(t),t,u,\lambda,\psi(t))\eta(t)(du)=\sup_{p\in
        P}\ct{H}(x_*(t),t,p,\lambda,\psi(t)).
   \end{equation}
   Рассмотрим множество ее всевозможных решений
   $(x,u,\lambda,\psi)\in\fr{X}\times\fr{U}\times\mm{R}\times\fr{X}$
    этой системы на $[0,n]$. Это множество компактно, поскольку выполнены
    условия продолжимости. Тогда компактно и множество всевозможных
    решений
   $(x,u,\psi,\psi_0)\in\fr{X}\times\fr{U}\times\fr{X}\times\mm{R}$
       обобщенной системы на $\mm{T}.$ Обозначим его через $\td{\fr{S}}.$

\begin{prop}
    В условиях $\bf{I}abc,\bf{II}$ если имеется
    оптимальное управление $u_*\in\fr{U}$  и соответствующая
    ему  $x_*=\ph(u_*)$, то  для некоторых $\Psi_0\in\mm{R},\psi\in\fr{X}$
    выполнено почти всюду на $\mm{T}$
    \rref{sys_x},\rref{sys_psi},\rref{maxH}.
%\begin{equation}
%   \label{sys}
%
%   \end{equation}
%
\end{prop}
\doc

% Заметим, что в силу липшицевости  $f,g$, равномерной на конечных промежутках
 %времени, правая часть включения ограничен, следовательно
 % всякое решение системы
 %\rref{sys},\rref{sys_psi},\rref{maxH} продолжимо на все $\mm{T}.$

 Рассмотрим для всякого $n\in\mm{N}$
 последовательность вспомогательных задач вида
 $$\dot{x}=f(t,x,u), u\in P, t\in [0,n] $$
 $$x(0)=0,\ x(n)=x^*(n),$$
  $$J_k=\int_{[0,n]} g(t,x(t),u(t))dt\to\max,$$

 Тогда $(x^*_n,u^*_n)$ --- оптимальна для вспомогательной задачи,
 отсюда  \cite[теорема 5.2.1]{clarke} для функции
$H$
 при некоторых $\lambda_n\in\mm{R},\psi_n\in C([0,n],\mm{R}^{m})$
 со свойством
 $\lambda_n+||\psi_n(0)||=1$
 при почти всех $t\in[0,n]$ выполнены соотношения \rref{sys_psi},
 \rref{maxH}. Продолжим $(x^*_n,u^*_n,\psi_n)$ с $[0,n]$ на $\mm{T}$
 как решение \rref{sys_x},\rref{sys_psi},\rref{maxH} произвольным
 образом.
% $$\dot{\psi}_n\in-\partial_x H(t,x^*,u^*,\psi_n, \psi^0_n)$$
% $$H(t,x^*,u^*,\psi_n, \psi^0_n)=\max_{p\in P} H(t,x^*,p,\psi_n, \psi^0_n).$$

 Заметим, дифференциальное включение \rref{sys_psi} полунепрерывно сверху по
 фазовой переменной, интегрально ограничено на ограниченных множествах, тогда
 графики решения этого включения не покидают за конечное время компактного
 множества, а в виду ограниченности - равностепенно непрерывны на
 всяком временном отрезке. Аналогично, $u^*_n|_{[0,n]}$ можно погрузить в
 $\Pi([0,n],P)$, а  $u^*_n$ в $\Pi(\mm{T},P).$
 Поскольку теперь кортежи $(x^*_n,u^*_n,\lambda_n,\psi_n)$
 погружены в компакт $\td{\fr{S}}$, то можно выделить подпоследовательность,
 сходящуюся к некоторому $(x^*,\eta,\Psi_0,\psi).$ Однако заметим,
 что в силу \rref{maxH} на промежутке $[0,n]$ точка $(g(t,x^*_n,u^*_n),f(x,u))$
 лежит на границе множества из условия ({\bf{III}}). Тогда в пределе
  точка $\int_P (g(t,x^*,u),f(x,u))\eta(t)(du)$ тоже
 лежит на границе множества из условия ({\bf{III}}), следовательно
 может быт реализована дираковской мерой, то есть
% в силу условия ({\bf{III}}) реализуется дираковской мерой,
 фактически управлением $u^*\in B(\mm{T},P)$.
 Но тогда $(x^*_n,u^*_n,\Psi^0_n,\psi_n)$
 сходятся к $(x^*,u^*,\Psi^0,\psi)$, и в силу замкнутости множества
 решений соотношений \rref{sys_x},\rref{sys_psi},\rref{maxH}
 эта четверка также им удовлетворяет.
  \bo

{\bf Условия трансверсальности}
  Cоотношения \rref{sys_x},\rref{sys_psi},\rref{maxH} не содержат условия на
  правом конце. Есть несколько вариантов таких условий (подробнее см.\cite[\S 1.6]{kr_as}),
  в данной работе исследуется
\begin{equation}
   \label{trans}
       \lim_{t\to\infty} \psi(t)=0.
   \end{equation}

%  Обозначим для каждого $u\in\fr{U}$ через $S(u)\subset\fr{X}$ множество всевозможных решений
%  линейного уравнения $\Psi$ при заданных $u,x=\ph{u}.$

% Условие $\bf{IV}$: Найдется такая окрестность $\ct{O}$
%  Существуют такие $\omega\in\Omega$ и окрестность
% траектории $x_{*},$ что для всякой траектории $x$ из неё,
%  для всякой удовлетворяющей \rref{sys},\rref{sys_psi},\rref{maxH} тройки
% $(u,\psi,\Psi_0)$
%  существует предел
%    $\lim_{t\to\infty} \psi(t)$ и для всех $s\in\mm{T}$
%    $\lim_{t\to\infty} ||\psi(s)-\psi(s)||\leq \omega(s)\omega(\rho(x,x^*)).$

 Условие $\bf{IV}$: Для всякой оптимальной для задачи {\bf{IV}} траектории $x^0$
 для всякого  решения $(x^0,u^0,\lambda^0,\psi^0)$ системы принципа максимума \rref{sys_x},\rref{sys_psi},\rref{maxH}
 найдется такая его окрестность $\Upsilon,$ в которой
 множитель Лагранжа $\psi^0$ устойчив, то есть
  для всякого $\epsi\in\mm{R}_{> 0}$ найдётся такое $\delta\in\mm{R}_{> 0}$  и
  $t\in\mm{T}$, что  если
  для решения  $(x,u,\lambda,\Psi)\in\Upsilon$ принципа максимума
\rref{sys_x},\rref{sys_psi},\rref{maxH}
выполнено
  $||\psi(t)-\psi^0(t)||<\delta$,
  $||x(t)-x^0(t)||<\delta$,
  $|\lambda-\lambda^0|<\delta$, то
  $||\psi(T)-\psi^0(T)||<\epsi$ имеет место для всех $T\in [t,\infty\rangle.$

%существует предел
%    $\lim_{t\to\infty} \psi^0(t),$ кроме того
% существуют такие $\omega\in\Omega$ и окрестность
% траектории $x_{*},$ что для всякой траектории $x$ из неё,
%  для всякой удовлетворяющей \rref{sys},\rref{sys_psi},\rref{maxH} тройки
%  $(u,\psi,\Psi_0)$
%  и для всех $s\in\mm{T}$
%    $\lim_{t\to\infty} ||\psi(t)-\psi(s)||\leq \omega(s).$

%   Заметим, что условия, гарантирующие $\bf{IV}$, имеются, например в
%   \cite[Lemma 4.1]{norv}, а основной результат работы \cite{norv}
%(в предположении единственности
%   оптимальной траектории)     является следствием

\begin{prop}
    В условиях $\bf{Iabc},\bf{II}-\bf{IV}$ для всякой неограниченно возрастающей
    последовательности
%   векторов $(l_n)_{n\in\mm{N}}\in (S^{m-1})^\mm{N}$ и
    моментов времени $(\tau_n)_{n\in\mm{N}}\in \mm{T}^\mm{N}$.
    найдется такая оптимальная для задачи \rref{opt} траектория $x^\infty$
    и такое решение
    $(x^\infty,u^\infty,\lambda^\infty,\Psi^\infty)$
 соотношений принципа максимума \rref{sys_x},\rref{sys_psi},\rref{maxH}, что   выполнено
\begin{equation}
   \label{partlim}
       \underline{\lim}_{n\to\infty}||\Psi^\infty(\tau_n)||_m=0.
   \end{equation}
%    найдётся такая тройка
%  $(u,\psi,\Psi_0)$ такое решение оптимальное решение $u_*\in\fr{U},$ что $x_*=\ph(u_*)\in\fr{X},$
 %   что
\end{prop}

\doc
  Зафиксируем некоторую неограниченно возрастающую
    последовательность
%   векторов $(l_n)_{n\in\mm{N}}\in (S^{m-1})^\mm{N}$ и
  моментов времени $(\tau_n)_{n\in\mm{N}}\in \mm{T}^\mm{N}$.

%  Для  всякого $n\in\mm{N}$ введем плоскость
%  $\alpha_n\rav\{x\in\mm{R}^m\,|\,l_n'x=l_n'x^0(\tau_n)\}.$

  Рассмотрим задачу
   $$J_{\tau_n}(u)\rav\int_{[0,\tau_n]} g(t,\ph(u)(t),u(t))dt\to\max.$$
  Она имеет решение $u^{n}\in\ct{U}$, обозначим порожденную этим управлением траекторию
   через $x^{n}\in\Phi.$

    Рассмотрим последовательность $(x^{n})_{n\in\mm{N}}.$ Она содержится в
  компакте $\Phi$,
  следовательно имеет хотя бы один частичный предел $x^{\infty}\in\Phi$
    в компактно-открытой топологии.
  Аналогично, прореживая последовательность управлений $u^{n}$, у
  нее можно
  также найти частичный предел  $u^{\infty}$, но уже в компакте {$\Pi(\mm{T}, P)$}.
  Более того, поскольку выполнено \rref{sys} при всяком $n\in\mm{N}$
  для пары $(x^{n},u^{n})$ на промежутке $[0,\tau_n]$, то для
  пределов  можно обеспечить
 $x^{\infty}=\td{\ph}(u^{\infty})$.

%  Поскольку по условию оптимальная траектория исходной задачи
%  единственна, то можно считать, что последовательность
%   $(x^{n})_{n\in\mm{N}}$ не покидает окрестности, указанной
%  в условии $\bf{IV}$ и  сходится к  $x^0=x^{\infty}=\td\ph(u^{\infty}).$

   Далее, поскольку все $x^{n}$ оптимальны в своих задачах,
    то для каждого на $[0,\tau_n]$ выполнен принцип максимума \rref{sys_x},\rref{sys_psi},\rref{maxH},
   причем с тем же самым гамильтонианом, что и для исходной задачи.
   В частности для некоторых  множителей Лагранжа
  $\lambda^{n}\in\mm{R},\Psi^{n}\in C([0,n],\mm{R}^m)$
   на $[0,\tau_n]$ имеют место соотношения принципа максимума
   при подстановке $x^*=x^{n},$ $u^*=u^{n}.$
%  со свойством $|\Psi^{n,l}_0|^2+||\Psi^{n,l}())||^2=1,$
%   так что выполнены основные соотношения для некоторого управления $u\in B([0,n],P).$
   Кроме того, имеет место условие трансверсальности на правом конце:
 $\Psi^{n}(\tau_n)=0.$

%  Для всякого $k\in\mm{N}$ последовательность $(\Psi^n|_{[0,k]})_{k\in\mm{N}}$
  Множество пар $(\lambda^n,\Psi^{n})$ содержится в компакте,
  следовательно имеет предельную точку. Очередной раз прореживая
  последовательность, можно считать, что она сходится.
  Поскольку соотношения принципа максимума полунепрерывно сверху
  зависят от $x,\psi,u$,  то этот предел $(\lambda^\infty,\psi^\infty)$
  вместе с оптимальными $x^\infty,u^\infty$   также удовлетворяет принципу максимума.
  В частности, в силу условия {\bf{III}} тогда можно считать, что $u^\infty\in \fr{U}.$

   Заметим, что $(u^{\infty},x^{\infty})$ оптимальна, действительно, исходная
  задача имеет некоторое оптимальное управление $u^0$, а в силу оптимальности $u^{n}$
  выполнено $J_{\tau_n}(u^{n})\geq J_{\tau_n}(u^0),$
  но левая часть неравенства  сходится к $\td{J}(u^{\infty})$, тогда как правая часть
  сходится к оптимальному решению исходной задачи $J(u^0).$

  Теперь в некоторой окрестности $(x^{\infty},u^{\infty},\lambda^\infty,\psi^\infty)$
   выполнено условие  {\bf{IV}}, в частности на нашей последовательности.
  Рассмотрим произвольное $\epsi\in\mm{R}_{> 0}$, согласно {\bf{IV}}
  возьмем $\delta$ и $t.$
  Поскольку найдется такой номер $N\in\mm{N},$
  что при $n>N$ выполнено $||\psi^n(t)-\psi^\infty(t)||<\delta,$
  то из {\bf{IV}} следует $||\psi^n(T)-\psi^\infty(T)||<\epsi$
   для всех $T\in [t,\infty\rangle.$
  В частности $||\psi^n(\tau_n)-\psi^\infty(\tau_n)||<\epsi$
  для всех $\tau_n\in [t,\infty\rangle,n>N$, то есть для всех $n$ начиная с
  некоторого номера $N'$.
  Но $\psi^n(\tau_n)=0$ для всех $n\in\mm{N}$, отсюда
  $||\psi^\infty(\tau_n)||<\epsi$ для всех $n>N'.$
  В силу произвольности $\epsi\in\mm{R}_{> 0}$
   \rref{partlim} показано.
%  Отсюда $\overline\lim_{n\to\infty} \psi^\infty(\tau_n)<\epsi$, то есть
%  $\lim_{n\to\infty} \psi^\infty(\tau_n)=0.$
\bo

Заметим, что условие {\bf{IV}} следует из такого условия:

Условие $\bf{IV}'$: Для всякой оптимальной для задачи \rref{opt}
траектории $x^0$
 для всякого  решения $(x^0,u^0,\lambda^0,\psi^0)$ системы принципа максимума
\rref{sys_x},\rref{sys_psi},\rref{maxH}
  существует предел
    $\lim_{t\to\infty} \psi^0(t),$ кроме того этот предел равномерен
    в целой окрестности для решений \rref{sys_x},\rref{sys_psi},\rref{maxH}, то есть
    найдутся такие функция $\omega_1\in\Omega$ и окрестность $\Upsilon$ точки
    $(x^0,u^0,\lambda^0,\psi^0)$, что
  для всякого решения  $(x,u,\lambda,\Psi)\in\Upsilon$
 принципа максимума \rref{sys_x},\rref{sys_psi},\rref{maxH}
   для всех $s\in\mm{T}$
    $\lim_{t\to\infty} ||\psi(t)-\psi(s)||\leq \omega_1(s).$

    (Действительно, в этих условиях у близкого к оптимальному решения автоматически
    на бесконечности имеется предел, в силу равномерности стремления к этому пределу
    $||\psi(T)-\psi^0(T)||$, $||\psi(\infty)-\psi^0(\infty)||$
    отличаются не более чем на $2\omega_1(T)$ на целой окрестности, осталось для
    $\delta=\omega_1(t)$ присвоить $\epsi=\max_{T>t}\omega_1(T)$).

\begin{sled}
    В условиях  $\bf{Iabc},\bf{II,III,IV'}$
    найдется такая оптимальная для задачи \rref{opt} траектория $x^\infty$
    и такое решение соотношений принципа максимума
    $(x^\infty,u^\infty,\lambda^\infty,\Psi^\infty)$, что   выполнено \rref{trans}.
\end{sled}
\doc
   Пусть не так, тогда поскольку предел     $\lim_{t\to\infty} \psi^\infty(t)$
   существует для всех интересующих нас $\psi^\infty$, то он отличен
   от нуля. Тогда, применяя предложение, получаем противоречие.
 \bo

%%%%%%%%%%%%%%%%%%%%%%%%%%%%%%%%%%%%%%%%%%%%%%%%%%%%%%%
%%%%%%%%%%%%%%%%%%%%%%%%%%%%%%%%%%%%%%%%%%%%%%%%%%%%%%%%%

%%%%%%%%%%%%%%%%%%%%%%%%%%%%%%%%%%%%%%%%%%%%%%%%%%%%%%%
%%%%%%%%%%%%%%%%%%%%%%%%%%%%%%%%%%%%%%%%%%%%%%%%%%%%%%%%%

%%%%%%%%%%%%%%%%%%%%%%%%%%%%%%%%%%%%%%%%%%%%%%%%%%%%%%%
%%%%%%%%%%%%%%%%%%%%%%%%%%%%%%%%%%%%%%%%%%%%%%%%%%%%%%%%%

%%%%%%%%%%%%%%%%%%%%%%%%%%%%%%%%%%%%%%%%%%%%%%%%%%%%%%%
%%%%%%%%%%%%%%%%%%%%%%%%%%%%%%%%%%%%%%%%%%%%%%%%%%%%%%%%%

 Условие $\bf{V}$: Для всякой оптимальной для задачи \rref{opt} траектории $x^0$
 для всякого  решения $(x^0,u^0,\lambda^0,\psi^0)$
  системы \rref{sys_x},\rref{sys_psi},\rref{maxH}
 найдется такая его окрестность $\Upsilon,$ в которой
 для всякого решения  $(x^0,u^0,\lambda,\Psi)\in\Upsilon$ соотношений \rref{sys_x},\rref{sys_psi}
 множитель Лагранжа $\psi$  устойчив для системы \rref{sys_x},\rref{sys_psi}.

\begin{prop}
    В условиях $\bf{Iabc},\bf{II}-\bf{III},\bf{V}$ для всякой неограниченно возрастающей
    последовательности
%   векторов $(l_n)_{n\in\mm{N}}\in (S^{m-1})^\mm{N}$ и
  моментов времени $(\tau_n)_{n\in\mm{N}}\in \mm{T}^\mm{N}$
    для всякой оптимальной для задачи \rref{opt} пары $(x^\infty,u^\infty)$
    найдется такое решение соотношений принципа максимума
    $(x^\infty,u^\infty,\lambda^\infty,\Psi^\infty)$, что   выполнено \rref{partlim}
\end{prop}

\doc

  Пространство $\Pi(\mm{T},P)$  является метризуемым,
  зафиксируем на нем некоторую метрику. Аналогично сделаем для $\mm{X}.$

  Каждому решению  $(x^\infty,u^\infty,\lambda^\infty,\psi^\infty)$
  системы \rref{sys_x},\rref{sys_psi},\rref{maxH}
  согласно {\bf{V}}
  можно сопоставить  свою окрестность $\Upsilon.$
  Всевозможные такие окрестности образуют покрытие $\fr{S}.$
  Поскольку последнее множество компактно, то
  из этого покрытия можно выделить конечное подпокрытие $\fr{S}'$, а
  расстояние от $\fr{S}$ до
  до границы $\fr{S}'$ будет не меньше некоторого $\varkappa\in\mm{R}_{>0}.$

  Зафиксируем произвольную оптимальную для задачи \rref{opt} пару $(x^\infty,u^\infty)$.
  Рассмотрим для каждого $n\in\mm{N}$ задачу
\begin{equation}
   \label{z_1}
   J_n(u)\rav\int_{[0,\infty\rangle}
   g(t,\ph(u)(t),u(t))-\frac{1}{n}e^{-t}||u(t)-u^{\infty}(t)|| dt
  \to\max.
   \end{equation}
  Ясно, что только пара $(x^\infty,u^\infty)$ является оптимальной для этой задачи.
  Обозначим через $\fr{S}_n$ множество четверок
  $(x^\infty,u^\infty,\lambda^n,\psi^n))$,
  удовлетворяющих
  \rref{sys_x},\rref{sys_psi} и
\begin{equation}
   \label{sys_max_}
        u^\infty(t)\in arg\max_{p\in
        P}
      \big(\ct{H}(x^n(t),t,p,\lambda,\psi(t))-\frac{1}{n}e^{-t}||p-u^{\infty}(t)|| \big).
   \end{equation}
 Заметим, что все $\fr{S}_n$ замкнуты, а поскольку содержатся в компакте, то и компактны.

Поскольку  соотношение \rref{sys_max_} также полунепрерывно
  сверху зависит от коэффициента перед последним слагаемым,
   все соотношения ограничены на ограниченных множествах,
  а множество из \rref{sys_max_} имеет сильный селектор, то в силу
  \cite[теорема 4.3.3]{tovst} имеет место
  для  пучков решений
   \rref{sys_x},\rref{sys_psi},\rref{sys_max_}
 полунепрерывность сверху  уже
   по коэффициенту перед последним слагаемым.
   В частности верхний предел компактов $\fr{S}_n$ вложен в $\fr{S}.$
   Следовательно, начиная с некоторого номера $N\in\mm{N}$, для всех $\fr{S}_n$
   будет являться покрытием конечное покрытие $\fr{S}'$.
   Тогда для каждого $n>N$ для новой задачи \rref{maxH_1}
   выполнено условие {\bf{IV}}.

   Рассмотрим некоторую неограниченно  возрастающую
    последовательность
%   векторов $(l_n)_{n\in\mm{N}}\in (S^{m-1})^\mm{N}$ и
  моментов времени $(\tau_k)_{k\in\mm{N}}\in \mm{T}^\mm{N}$.

 Теперь по предыдущему предложению для каждого  $n\in\mm{N}$
 при некотором $(x^\infty,u^\infty,\lambda^n,\psi^n)\in\fr{S}_n$
 выполнено $\underline{\lim}_{k\to\infty}\psi^n(\tau_k)=0.$
 Кроме того, по уже показанному последовательность
 $(x^\infty,u^\infty,\lambda^n,\psi^n)_{n\in\mm{N}}$
 имеет предельную точку
  $(x^\infty,u^\infty,\lambda^n,\psi^n)\in \fr{S}.$

   Рассмотрим произвольное $\epsi\in\mm{R}_{> 0}$, по условию {\bf{V}}
  найдутся такие  $\delta\in\mm{R}_{> 0}$ и $t\in\mm{T}$, что из
  $||\psi^n(t)-\psi^\infty(t)||<\delta$
  следует $||\psi^n(T)-\psi^\infty(T)||<\epsi$
   для всех $T\in [t,\infty\rangle.$
   Поскольку для $(\psi^n)_{n\in\mm{N}}$ точка $\psi^\infty$ ---
   предельная, то для некоторого $N\in\mm{N}$
  $||\psi^N(t)-\psi^\infty(t)||<\delta$, то есть
  $||\psi^N(T)-\psi^\infty(T)||<\epsi$. Тогда
  $\underline{\lim}_{k\to\infty}\psi^\infty(\tau_k)<\epsi.$
  В силу произвольности $\epsi\in\mm{R}_{> 0}$
  показано \rref{partlim}. \bo

 Аналогично, условие {\bf{V}} следует из такого условия:

Условие $\bf{V}'$: Для всякой оптимальной для задачи \rref{opt}
траектории $x^0$
 для всякого  решения $(x^0,u^0,\lambda^0,\psi^0)$ системы принципа максимума
\rref{sys_x},\rref{sys_psi},\rref{maxH}
  существует предел
    $\lim_{t\to\infty} \psi^0(t),$ кроме того этот предел равномерен
    в целой окрестности для решений \rref{sys_x},\rref{sys_psi}, то есть
    найдутся такие функция $\omega_1\in\Omega$ и окрестность $\Upsilon$ точки
    $(x^0,u^0,\lambda^0,\psi^0)$, что
  для всякого решения  $(x,u,\lambda,\Psi)\in\Upsilon$
 принципа максимума \rref{sys_x},\rref{sys_psi}
   для всех $s\in\mm{T}$
    $\lim_{t\to\infty} ||\psi(t)-\psi(s)||\leq \omega_1(s).$

\begin{sled}
    В условиях  $\bf{Iabc},\bf{II,III,V'}$
    найдется такая оптимальная для задачи \rref{opt} траектория $x^\infty$
    и такое решение соотношений принципа максимума
    $(x^\infty,u^\infty,\lambda^\infty,\Psi^\infty)$, что   выполнено \rref{trans}.
\end{sled}

 Одно из самых общих условий на \rref{trans} показано в  \cite[Theorem 6.1]{norv}.
 Если ограничиться лишь задачей управления без фазовых ограничений,
 то, поскольку
 условие {\bf{IV'}}
 следует   из  условий \cite[Theorem 3.1]{norv} в силу
  \cite[Lemm 3.1]{norv}, то и  основной результат статьи \cite{norv}
  вкладывается в следствие 3.

  Заметим, что при быстро убывающей функции $\omega$ из свойства {\bf{II}}
  многие множители Лагранжа будут удовлетворять \rref{trans} и само
  по себе это условие не позволит выделить существенно меньшее
  семейство экстремалей. Однако, как замечено в
  \cite[Theorem 8.1]{norv}, это можно исправить
  (\cite[Example  10.2]{norv}) усилив
  свойство  \rref{trans}.

  Заметим, что в условиях предложений 2,3 показано большее, что
  удовлетворяющий
  \rref{partlim} множитель Лагранжа $\Psi$ является пределом отображений
  $\Psi_n$, зануляющихся на все больших моментах времени.
  В частности
\begin{sled}
    В условиях $\bf{I}-\bf{III},\bf{IV'}$
    найдется такая оптимальная для задачи \rref{opt} пара $(u^\infty,x^\infty)$
    и такое решение соотношений принципа максимума
    $(x^\infty,u^\infty,\lambda^\infty,\Psi^\infty)$, что   выполнено
    \rref{trans} и для некоторой возрастающей
    последовательности $(T_n)_{n\in\mm{N}}\in\mm{T}^\mm{N}$
    имеется сходящаяся к $(x^\infty,u^\infty,\lambda^\infty,\Psi^\infty)$
    последовательность $(x^n,u^n,\lambda^n,\Psi^n)_{n\in\mm{N}}$
    решений \rref{sys_x},\rref{sys_psi},\rref{maxH}
    со свойством $\Psi^n(T_n)=0.$
\end{sled}
\begin{sled}
    В условиях $\bf{I}-\bf{III},\bf{V'}$
    для всякой оптимальной для задачи \rref{opt} пары $(u^\infty,x^\infty)$
    найдется такое решение соотношений принципа максимума
    $(x^\infty,u^\infty,\lambda^\infty,\Psi^\infty)$, что   выполнено
    \rref{trans} и для некоторой возрастающей
    последовательности $(T_n)_{n\in\mm{N}}\in\mm{T}^\mm{N}$
    имеется сходящаяся к $(x^\infty,u^\infty,\lambda^\infty,\Psi^\infty)$
    последовательность $(x^n,u^n,\lambda^n,\Psi^n)_{n\in\mm{N}}\in\fr{S}^\mm{N}$
    решений \rref{sys_x},\rref{sys_psi}
    со свойством $\Psi^n(T_n)=0.$
\end{sled}
%
%возможно для некоторых линейных систем тогда автоматически выпишется формула
% то есть пишем решение краевой задачи и переходим к пределу
%
\begin{zam}
   Заметим, что если  в условиях двух последних следствий
   потребовать также равномерную ограниченность $x$ на оптимальных
    траекториях, то кроме     \rref{trans}
   автоматически будет выполнено также более тонкое условие
   транверсальности
   $$\lim_{t\to\infty} \Psi^\infty(t)'x^\infty(t)=0.$$
\end{zam}

%  Для задачи управления без фазовых ограничений показано чуть больше, чем в
%  \cite[Theorem 8.1]{norv} (сходимость всей четверки в
%   компактно открытой топологии, а не сходимость лишь $\Psi$ и поточечно).

% Возможно условия этого следствия близки к достаточным
% ведь существование такой последовательности автоматически
% означает оптимальность предела
% интересно лишь выяснить, можно ли всякую оптимальную траекторию
% (а не какую-то из них) представить в таком виде

 Заметим также, что если  для элементов так построенной
 последовательности $(x^n,u^n,\lambda^n,\Psi^n)_{n\in\mm{N}}\in\fr{S}^\mm{N}$
 выполнены достаточные условия оптимальности, то всякая ее предельная
 точка является оптимальным решением \rref{opt}.

  Для распространения результатов типа предложения 2 на случай задач
  с фазовыми ограничениями основной сложностью по-видимому
  будет доказательство того, что пучок условий принципа максимума
  замкнут.

%  (В )

%можно привести пример системыдлякоторой trans не будет partlim будет
% множители ведут себя как пружины.
% пружины отличаются параметром, он отвечает за сдвиг (точнее угол сдвига) "оси"
% пружины и фазу пружины.
% Ось пружины ведет асимптотически идет к нулю
% каждая пружина касается ноля в счетном числе точек
% разные фазы гарантируют невыполнение trans
%
%

\end{document}